\documentclass[11pt]{article}
\usepackage{amssymb}
\usepackage{latexsym}
\newtheorem{thm}{Theorem}[section]
\newtheorem{prop}[thm]{Proposition} \newtheorem{lemma}[thm]{Lemma}
\newtheorem{cor}[thm]{Corollary} \newtheorem{dfn}[thm]{Definition}
 \newtheorem{rmk}[thm]{Remark}
\newtheorem{ex}[thm]{Example} 

\newcommand {\pf}{\noindent{\bf Proof.}\ }

\newcommand{\naturals}{{\mathbb N}}
\newcommand{\reals}{{\mathbb R}}

\newcommand{\integers}{{\mathbb Z}}

\newcommand{\cala}{{\cal A}}

\newcommand{\calo}{{\cal O}}

\newcommand{\calv}{{\cal V}}
\newcommand{\calu}{{\cal U}}
\newcommand{\calw}{{\cal W}}

\newcommand{\qed}{\begin{flushright} $\Box$\ \ \ \ \ \end{flushright}}

\newcommand{\frakb}{\mathfrak{b}}
\newcommand{\frakg}{\mathfrak{g}}

\newcommand{\arrows}{\,\lower1pt\hbox{$\longrightarrow$}\hskip-.24in\raise2pt
             \hbox{$\longrightarrow$}\,}
\title{{\bf Linearization of Regular Proper Groupoids}}
\author{Alan
Weinstein\thanks{Research partially supported by NSF Grants
DMS-96-25122 and DMS-99-71505
and the Miller Institute for Basic Research in Science.
\newline \mbox{~~~~}MSC2000 Subject Classification Numbers: 58H05 (Primary), 57R99
             (Secondary).
\newline \mbox{~~~~}Keywords: Lie groupoid, proper action, submersion.
}
\\Department of Mathematics\\ University of California\\ Berkeley, CA
94720 USA\\ {\small(alanw@math.berkeley.edu)}}
\date{June 25, 2001}
\begin{document}
\maketitle
\begin{abstract}
Let $G$ be a Lie groupoid over $M$ such that the target-source map
from $G$ to $M \times M$ is proper.  We show that, if ${\mathcal O}$ is
an orbit of finite type (i.e. which admits a proper function with finitely many critical
points), then the restriction $G|_{\mathcal U}$
 of $G$ to some
neighborhood ${\mathcal U}$ of ${\mathcal
O}$ in $M$ is isomorphic to a similar restriction of the action groupoid for
the linear action of the transitive groupoid $G|_{\mathcal O}$
on the normal bundle $N\calo$.  The proof uses a deformation argument
based on a cohomology vanishing theorem, along with a slice
theorem which is derived from a new result on submersions with a fibre
of finite type.
\end{abstract}

\section{Introduction}
\label{sec-intro}
A Lie groupoid $G\arrows X$ is called {\bf proper} if the
(target,source) map $G\to X\times X$ is a proper map.  Such groupoids
arise, for example, as the transformation groupoids attached to smooth
proper actions of groups.  For proper group actions, and hence for
these transformation groupoids, there is a normal form valid in the
neighborhood of each orbit.    (We refer to the first two
chapters of \cite{du-ko:lie} for a detailed treatment of the theory of
smooth proper actions.)  Transformation to the normal form depends on: 
(1) the existence
of slices, which relate the behavior near an orbit to that near a
fixed point of the isotropy subgroup; (2) linearizability of actions
of compact groups near fixed points.  Together, these two ingredients
show that a proper group action is equivalent
in a neighborhood of each orbit to a linear action on a vector bundle.

As is explained in \cite{we:linearization}, we hope eventually to
establish a linearization theorem for all proper groupoids near
 their orbits.  Its most
general version would require a still unproven linearization
theorem around fixed points,
so in the present paper we sidestep this
problem by assuming that the groupoid is {\bf regular} in the sense
that its orbits all have the same dimension.  

The assumption that $G$ is regular greatly simplifies the local
problem, since the restriction $G_{\Sigma}$ of a
regular groupoid $G$ to a slice $\Sigma$ is {\bf essentially {\'e}tale}
in the sense that its action on the slice factors through an {\'e}tale
groupoid.  The {\'e}tale groupoid is easily linearized near its fixed
point, and a deformation argument using Crainic's cohomology vanishing
theorem \cite{cr:differentiable} allows us to linearize the
essentially \'etale groupoid
 $G_{\Sigma}$;  i.e., we prove that $G_{\Sigma}$ is
locally isomorphic to the action groupoid for the linear action of the
isotropy group $G_x$ on the tangent space $T_x\Sigma.$  

The second building block of our normal form is the restriction
of $G$ to an orbit, which is transitive and
is therefore the gauge groupoid of a principal bundle.  Combining the
action groupoid over a slice with the gauge groupoid over an orbit
produces the action groupoid for the (linear) action of the gauge
groupoid on the normal bundle of the orbit, which we call the
linear approximation to $G$ along the orbit.  Our main theorem asserts
that, under a differential-topological finiteness
assumption which we describe in the
following paragraph, $G$
is isomorphic to its linear approximation in some neighborhood of each
orbit.  We refer to such an isomorphism as a linearization of $G$ along
the orbit.

The main theorem follows from a slice theorem which asserts that a
linearization of the restriction to a slice extends to a linearization
along an orbit.  This slice theorem holds even in the nonregular
case, but it turns out to require a
differential-topological assumption on the orbit itself. The necessity
of this assumption appears already in the special case where the proper
groupoid has trivial isotropy groups, i.e. where $G$ is the
equivalence relation whose
equivalence classes are the fibres of a submersion $f:X \to Y$.
Applied to such  a groupoid, our normal form theorem asserts that,
restricted to some open neighborhood $\calu$ of each fibre
$f^{-1}(y)$, $f$ is a trivial fibration onto $f(\calu)$.  Already in
this case, it turns out that an additional hypothesis is necessary;
the most natural one seems to be that the fibre $f^{-1}(y)$ is of
{\bf finite type} in the sense that it admits a proper map to $\reals$ with
finitely many critical points.  (See Appendix \ref{sec-finite} for
a discussion of this finiteness condition.) 

The body of the paper begins with a discussion of the definition of  proper
groupoids, followed by some examples showing the importance of a local
triviality condition.  We prove our linearization theorem for
groupoids associated with submersions, and then for {\'e}tale
groupoids.  After proving a rigidity theorem for proper groupoids, we
deal with the effectively {\'e}tale case.  Finally, we prove the slice
theorem and use it to deduce the main theorem.
Two appendices are devoted to background material on proper
mappings and manifolds of finite type.

Proper groupoids seem to 
have appeared only infrequently in the literature. Moerdijk and Pronk
\cite{mo-pr:orbifolds} characterized orbifolds as Morita
equivalence classes of {\'e}tale proper groupoids on manifolds.  The
orbit spaces of regular proper groupoids are orbifolds as well, and
the Morita equivalence classes of these groupoids may be thought of as
principal bundles over orbifolds.  In connection with operator
algebras, proper groupoids
 appear in Connes' book \cite{co:noncommutative} as the
appropriate setting for the geometric realization of cycles in
$K$-theory; an appendix (Section 6) on proper groupoids in the related 
paper of Tu \cite{tu:hyperboliques}
includes the construction of a cutoff function used
by Crainic \cite{cr:differentiable} in the proof of his cohomology
vanishing theorem.   

I would like to thank Ronnie Brown, Marius Crainic, Bob Gompf, Rob
Kirby, Kirill Mackenzie, Ga{\"e}l Meigniez, and Ieke Moerdijk for
their helpful advice concerning groupoids on the one hand and
manifolds of finite type on the other.  Comments 
by Moerdijk and the referee on the first version of this paper were
particularly useful.  For
hospitality during the final stages of this work, I would like to
thank Harold Rosenberg at l'Institut de Math{\'e}matiques de Jussieu,
Yvette Kosmann-Schwarzbach and Claude Viterbo at l'{\'E}cole
Polytechnique, and Peter Michor at the Erwin Schr\"odinger Institute.

\section{Definition of proper groupoids}
\label{sec-definition}

For basic notions about groupoids, we refer the reader to
\cite{ca-we:geometric} or \cite{ma:lie}.
Our conventions and notation include the following.
``Lie groupoid'' will always mean ``smooth groupoid,'' not necessarily
transitive as was the case in \cite{ma:lie}.
  All manifolds will be Hausdorff;
Appendix \ref{sec-propermapping} shows why this assumption is
reasonable in the context of proper groupoids.  All neighborhoods will
be open.  In a groupoid $G\arrows X$, we will denote the target
and source maps by $\alpha$ and $\beta$ respectively.  We will
sometimes call the map $(\alpha,\beta):G\to X\times X$ the {\bf anchor} of the
groupoid, following \cite{ma:lie}.  If $A$ and $B$ are subsets of $X$,
we denote by $G_{AB}$ the set $\alpha^{-1}(A)\cap\beta^{-1}(B)$.  If
$A=B$, we write $G_A$ for $G_{AA}$.  For one-point subsets, we
abbreviate $G_{\{x\}\{y\}}$ by $G_{xy}$.  Combining the two
abbreviations leads to the usual notation $G_x$ for the isotropy group
of $x\in X$.

A submanifold $S\subseteq G$ for which $(\alpha,\beta)(S)$ is the
graph of a diffeomorphism $\phi_S:X\to X$ is called a {\bf bisection} of
$G$.  The bisections form a group under setwise $G$-multiplication
which acts on $G$ by left or right translations and on $X$ via the
homomorphism $S\mapsto \phi_S$.

Following \cite{mo-mr:integrability},
we will say that a groupoid $G\arrows M$ is ``source-xxx'' if the
target and source maps of $G$ (or their fibres, as will be clear from
the context) each have the property ``xxx.''  For instance, we may
refer to groupoids as being source-proper, source-connected, etc.

A groupoid $G$ is {\bf regular} if its orbits all have the same
dimension.  This condition is equivalent to constancy of rank for
either the anchor of $G$ or the anchor $\cala(G)\to TX$ of its
Lie algebroid.

We now recall two definitions concerning group actions.

\begin{dfn}
\label{dfn-properaction}
An action of a Lie group $\Gamma$ on a topological space $X$ is {\bf
proper} if the mapping $(\gamma,x)\mapsto (\gamma x,x)$ is a proper
mapping from $\Gamma \times X$ to $X\times X$.
\end{dfn}

\begin{dfn}
\label{dfn-actiongroupoid}
Given an action of a group $\Gamma $ with identity element $e$ on a set $X$,
the corresponding {\bf action groupoid} is the groupoid $\Gamma\times
X\arrows X$ with target and source maps
$(\gamma,x)\mapsto \gamma x$ and $(\gamma,x) \mapsto x$, product
$(\gamma,\gamma'x)(\gamma',x)=(\gamma \gamma',x)$, unit embedding
$x\mapsto (e,x)$, and 
inversion $(\gamma,x)^{-1}=(\gamma^{-1},\gamma x)$.
\end{dfn}

On the basis of these two definitions, we introduce the notion of
proper groupoid.  

\begin{dfn}
\label{dfn-propergroupoid}
A {\bf proper groupoid} is a Lie groupoid $G\arrows X$ for
which the anchor mapping
$(\alpha,\beta):G\to X\times X$ is proper. 
\end{dfn}

Here is an important property of proper groupoids.  

\begin{prop}
\label{prop-closed}
Each orbit of a proper groupoid is a closed submanifold.
\end{prop}
\pf
Let $\calo$ be the $G$-orbit through $x\in X$.  The isotropy $G_x$
is a compact group acting freely on $G_{xX}$ by multiplication
from the left, and $\beta$ factors through the natural projection
to give a map from the quotient $G_{xX}/G_x$ to $X$ which is an
injective immersion with image $\calo$.  To show that $\calo$ is
closed, it suffices to show that this immersion is proper.  

To this end, let $g_i$ be a sequence of elements in $G_{xX}$
such that $\beta{g_i}$ is convergent.  Then the anchor of $G$ applied to $g_i$ 
gives the convergent sequence $(x,\beta{g_i})$.  Since the groupoid is proper, $g_i$
contains a convergent sequence, hence so does the corresponding sequence $[g_i]$ in
the quotient space   $G_{xX}/G_x$.
\qed

It turns out that properness of a groupoid is not sufficient to imply
some of the nice properties which we associate with proper actions of
groups.  For this reason, we will sometimes impose the additional
condition of {\bf source-local triviality}.
Examples in Section \ref{sec-stability} below will show that this
condition does not follow from the properness of
the anchor.  For action groupoids associated with group actions,
the target and source maps are globally, hence locally, trivial
fibrations.

\section{Groupoid actions and stability}
\label{sec-stability}

Let $G\arrows X$ be a groupoid, and let $\mu:Y\to X$ be a surjective
mapping.  An {\bf action} of $G$ on $Y$ is a mapping
$(g,y)\mapsto gy$ to $Y$ from the fibre product $G\times_X Y$ (using
the source map from $G$ to $X$) satisfying the usual conditions for
associativity and action of the identities.  

\begin{ex}
\label{ex-action}
{\em
Any groupoid $G\arrows X$ acts on its base
$X$ by the rule $gx=\alpha(g)$ whenever $\beta(g)=x$.  If $H$ and $A$
are subsets of $G$ and $X$ respectively, $HA$ is defined in the usual
way, as in the case of group actions.  
}
\end{ex}

Definition \ref{dfn-actiongroupoid} is easily extended from group
actions to groupoid actions.
Given an action of $G\arrows X$ on $Y$,
the associated {\bf action groupoid} is the groupoid
$G\times_X Y\arrows Y$ with anchor $(g,y)\mapsto (gy,y)$ and
multiplication $(h,gy)(g,y)=(hg,y)$.  The mapping $\mu$ is sometimes
called the {\bf moment map} of the groupoid action.  In the 
differentiable category, the fibre product is a manifold because the
source map of $G$ is a submersion, so that
we may require the action to be differentiable, in which case the
action groupoid $G\times_X Y\arrows Y$ is again a Lie groupoid.
The action groupoid for the action of $G$ on its base $X$ is
naturally isomorphic to $G$ itself.

Motivated by the case of group actions, we make the following
definition.
\begin{dfn}
\label{dfn-stable}
A fixed point $x$ of a topological groupoid $G\arrows X$ is {\bf
stable} if every neighborhood of $x$ contains a $G$-invariant
neighborhood.
\end{dfn}

The following theorem shows that proper groupoids share an important
property with proper actions of groups.

\begin{thm}
\label{thm-stability}
Every fixed point of a source-locally trivial proper topological
groupoid is stable.
\end{thm}
\pf
Let $x$ be a fixed point of $G\arrows X$.  To show that a given
neighborhood $\calu$ of $x$ contains an invariant neighborhood, we may
assume to begin that $\calu$ is small enough so that
$\beta^{-1}(\calu) \approx \calu\times G_x$ as spaces over $\calu$.
Since $G_x$ is compact, each neighborhood of $\beta^{-1}(y)$ for $y$
in $\calu$ contains a neighborhood of the form $\beta^{-1}(\calv)$.  

Now we define the {\bf core} $C(\calu)$ of $\calu$ to be 
$$\{y\in \calu | Gy \subseteq \calu\} = \{y\in
\calu|\beta^{-1}(y)\subseteq \alpha^{-1}(\calu)\}.$$
$C(\calu)$ is invariant: if $y\in C(\calu)$ and $gy $ is defined,
then $h(gy)=(hg)y\in \calu$ whenever $h(gy)$ is defined, so $gy\in C(\calu)$.
$C(\calu)$ is open: if $y\in C(\calu)$, then $\beta^{-1}(y)$ is
contained in the open set $\alpha^{-1}(\calu)$, so $y$ has a
neighborhood $\calv$ such that $\beta^{-1}(\calv)\subseteq
\alpha^{-1}(\calu)$, i.e. $\calv\subseteq C(\calu)$.
\qed

Examples \ref{ex-twocharts} and \ref{ex-dipole} below show that
Theorem \ref{thm-stability} requires an assumption like source-local
triviality.  These examples also show that, while the condition of
properness is preserved under Morita equivalence of groupoids
\cite{mo-mr:notes},
source-local triviality and stability are not.

\begin{ex}
\label{ex-twocharts}
{\em Let $X_1$ be the plane $\reals^2$ with the origin removed.
Let the groupoid
$G_1$ be the equivalence relation on $X_1$, with quotient space
$\reals$, consisting of all the pairs of points lying on the same
vertical line.  The anchor of $G_1$ is proper, but the source map is
not locally trivial over any point on the vertical line through the
origin.  There are no fixed points.

Now let $G_2$ be the restriction of $G_1$ to any
horizontal line $X_2$ {\em not} passing through
the origin.  $G_2$ is just the trivial groupoid, with all elements
units, so it is source-locally trivial as well as proper.
Every element of $X_2$ is a stable fixed point for $G_2$.
The two groupoids are Morita equivalent since the second is
the restriction of the first to
a closed submanifold passing through all orbits and intersecting them
transversely. (See  Example 2.7 in \cite{mu-re-wi:equivalence} for the
relevant topological result; the smooth case is handled similarly).

Finally, let $G_3$ be the restriction of $G_1$ to the union $X_3$ of {\em two}
horizontal lines, one of which passes through (but does not include)
the origin.  This
groupoid is equivalent to the first two; like $G_1$ it has a proper
anchor map, but it is not source-locally trivial.  The point of
$X_3$ lying on the vertical line through the origin is a fixed point
for $G_3$ (the only one), but it is not stable.  
}
\end{ex}

\begin{ex}[Dipole foliation]
\label{ex-dipole}
{\em
Let $X$ be the plane with the two points $(0,1)$ and $(0,-1)$ removed.  In
this plane, consider the  foliation given by the level curves of the
potential function produced by a unit positive charge at one deleted point
and a unit negative charge at the other.  This foliation is symmetric
about the horizontal axis; its
leaves are this axis and simple closed curves surrounding the two
deleted points.  Let $G$ 
the equivalence relation determined by the foliation, considered as a
groupoid over $X$.  
$G$ is easily seen to be proper; for instance, it is Morita
equivalent to its restriction to an open line segment joining the two
charges, which is a trivial groupoid.  On the other hand, $G$ is not
source-locally trivial around points of the horizontal axis.  
  We do obtain an equivalent source-locally trivial groupoid
by restricting $G$ to the open strip between the horizontal lines
through the charges.  The latter groupoid is in fact isomorphic to an
action groupoid for an $\reals$ action.
}
\end{ex}

\section{The main theorem}
\label{sec-main}

In this section, we state our main theorem.  Succeeding sections are
devoted to proofs of special cases, culminating in the proof of the
theorem itself.

For any Lie groupoid  $G\arrows X$, its 1-jet prolongation acts
naturally on $TX$, and this action descends to an action of $G$ itself
on the normal ``bundle'' to the
singular foliation of $X$ by orbits.
(See, for example, Appendix B in
\cite{ev-lu-we:transverse}.)  Consequently, for any orbit
$\calo$, the restricted groupoid $G_{\calo}$ has a
natural representation on the normal bundle $N\calo$.  The action
groupoid $G_{\calo}\times_{\calo} N\calo$ should be thought of as the
linear approximation to $G$ along $\calo$, so the following theorem
states that, near an orbit of finite type, a proper groupoid is
isomorphic to its linear approximation.

\begin{thm}
\label{thm-main}
Let $G\arrows X$ be a regular, proper Lie groupoid, and let $\calo$ be
an orbit of $G$ which is a manifold
of finite type.  Then there is a neighborhood
$\calu$ of $\calo$ in 
$X$ such that the restriction of $G$ to $\calu$ is isomorphic to the
restriction of the action groupoid $G_{\calo}\times_{\calo}
N\calo$ to a neighborhood of the zero section in $N\calo$. 
\end{thm}

\section{Local semitriviality of submersions}
\label{sec-submersions}

The prefix ``semi'' in the title of this section is meant in the same
sense as in ``semicontinuous''.  

If $f:X\to Y$ is a submersion, then
the fibre product $X\times_Y X$ is an equivalence relation, as well as
being a submanifold of $X\times X$.  With these structures, $X\times_Y
X$ becomes a proper Lie subgroupoid of $X\times X$.  (In fact, these
are the only proper Lie subgroupoids of $X\times X$.) The isotropy
groups of $X\times_Y X$ are trivial and the orbits are the fibres of
$f$. In this situation, our main theorem is essentially 
equivalent to Theorem \ref{thm-submersion} below, of which Example
\ref{ex-dipole} is an illustration.
The equivariant case is included for
later use.

\begin{thm}
\label{thm-submersion}
Let $f:X\to Y$ be a submersion.  For any $y\in Y$, if
$\calo=f^{-1}(y)$ is a manifold of
finite type, then there is a
neighborhood $\calu$ of $\calo$ in $X$ such that
$f|_{\calu}:\calu\to f(\calu)$ is a trivial fibration.  In other
words, there is a retraction $\rho:\calu\to \calo $ such that
$(\rho,f):\calu\to \calo \times f(\calu) $ is a diffeomorphism.

If $f$ is equivariant with respect to actions of a compact group $K$
on $X$ and $Y$, with $y$ a fixed point, and $f^{-1}(y)$ is of finite type
as a $K$-manifold, then $\calu$ can be chosen to be $K$-invariant and
$\rho$ to be $K$-equivariant.
\end{thm}
\pf
All which follows can be
done equivariantly.  
Let $h:\calo \to [0,\infty)$ be a proper function whose
critical points form a compact set.  Choose $N$ large enough so that
all the 
critical points of $h$ lie inside $h^{-1}([0,N-1))$.
Since $\calo $ is a closed submanifold of $X$, it has a
tubular neighborhood on which there is a smooth retraction 
$\rho_0$ to $\calo$.  Since the problem is local around $\calo$, we
will assume from now on that this neighborhood is all of $X$.
By composition with $\rho_0$, we extend
the function $h$ to $X$ and call the extension $h$ as well. 

For any $N>0$ and any subset $A$ of $X$, we will denote the 
subset of $A$ on which $h(x)<N$ by $A_{(N)}$ and the subset on
which $h(x)\leq N$ by $A_{[N]}$. 

Since $\calo_{[N+1]}$ is compact, we can choose a small open disc
$\calw$ around $y$ in $Y$ such that the restriction of $(f,\rho_0)$ to
$f^{-1}(\calw)_{[N+1]}$ is a diffeomorphism from
$f^{-1}(\calw)_{[N+1]}$ to $\calw\times \calo_{[N+1]}$.  We denote the
inverse of this diffeomorphism by $\Psi_0$.
By choosing $\calw$ small enough, we can also assure that the
function $h$ has no fibrewise 
(for $f$) critical points in $f^{-1}(\calw)_{[N+1]}\setminus
f^{-1}(\calw)_{(N-1)}$.

Now let $\calv$ be the union of $f^{-1}(\calw)_{(N+1)}$ 
and the set of points in $ f^{-1}(\calw)$ where
$h$ is not critical along the fibres of $f$.  It is a neighborhood of
$\calo$.  Choose (see Lemma
\ref{lemma-complete}) a complete riemannian metric on $\calv$.  With
the help of a compactly supported function equal to 1 on
$f^{-1}(\calw)_{[N+1]}$, we may modify this metric without affecting
its completeness so that its restriction to $f^{-1}(\calw)_{[N+1]}$
becomes a product metric with respect to the diffeomorphism with
$\calw \times \calo_{[N+1]}$.  As a result, the projection $\rho_0$
becomes a local isometry to $\calo$ when further restricted to each
fibre of $f$.

Let
$\xi$ be the vector field on $\calv$ which is the fibrewise gradient
of $h$.  Since the metric restricted to each fibre is complete, and
$h$ is proper, $\xi$ is a complete vector field.  Denote the flow of
$\xi$ by $t\mapsto\phi_t$, and let $\calu$ be the union of all the
trajectories of $\xi$ which intersect $f^{-1}(\calw)_{(N+1)}.$
  Since $\calu$ is the union over all $t$ of
$\phi_t(f^{-1}(\calw)_{(N+1)})$, it is an open subset of $X$.  It
contains $\calo$, and, since $h$ has no critical values in $\calu$
above $N-1$,
every point $x$ of $\calu$ for which $h(x)>N$
flows down under the flow of $-\xi$  to the level $h^{-1}(N)$. 

Now we define a retraction $\rho:\calu \to \calo$ as follows.  If $h(x)
< N+1$, we let $\rho(x)=\rho_0(x)$, where $\rho_0(x)$ is the original
tubular neighborhood retraction on $X$.  If $h(x) > N$, there is a
unique positive number $\tau(x)$, depending smoothly on $x$, such that
$h(\phi_{ -\tau (x)}(x)) = N$.  We then define $\rho (x) =
\phi_{\tau(x)}(\rho_0 (\phi_{ -\tau (x)}(x))).$
By the product structure of the metric, the map
$\rho_0$ commutes with the flow on $f^{-1}(\calw)_{[N+1]}$, 
and so these two definitions agree on
the intersection of the domains.  It is also clear that $\rho$ is a
retraction.  

To see that the product map $(\rho,f)$ is a
diffeomorphism from $\calu$ to $\calw \times \calo$, we observe that
it has an inverse $\Psi$ defined on $\calw\times \calo_{(N+1)}$ by
$\Psi=\Psi_0$ and on $\calw\times (\calo\setminus \calo_{[N]})$ by
$\Psi(y,z) = \phi_{\tau(z)}(\Psi_0(y,\phi_{-\tau(z)}(z))).$    
\qed

A proper submersion is always a locally trivial fibration.  
Compactness of all the fibres is not enough to insure
local triviality (restrict Example \ref{ex-tracks} below to $(0,1)\times
\naturals$), but it is an easy corollary of 
Theorem \ref{thm-submersion} that compactness
and connectedness of the fibres {\em are} enough.  These results were
already essentially proved by Ehresmann \cite{eh:espaces}. 
Meigniez \cite{me:submersions1} \cite{me:submersions2} 
has recently given several sufficient
conditions for submersions with noncompact fibres to be locally trivial.

Here is a very simple example which shows that the hypothesis
of finite type cannot be omitted from Theorem \ref{thm-submersion}.

\begin{ex}[Railroad to infinity]\label{ex-tracks}
{\em
Let $X \subset (-1,1)\times \naturals$ be the subset $\{(x,n)|nx^2<1\}$,
$Y=(-1,1),$ and $f$ the projection on the first factor.  Then there is
no neighborhood of $f^{-1}(0)$ on which $f$ becomes a trivial
fibration, since the inverse image of any $y\neq 0$ is finite and
hence cannot contain an embedded copy of $f^{-1}(0)$.
}
\end{ex}

Lest the reader think that connectedness might help, we give another example.

\begin{ex}[Ladder to heaven]\label{ex-ladder}
{\em Let $S$ be an infinite-holed torus with one end, and let $h:S\to [0,\infty)$ 
be a
proper function
for which the inverse image of each interval $[0,n]$ is a surface of
finite genus bounded by two circles.  Now let $X\subset (-1,1) \times
S$ be the subset $\{(x,s)|h(s)x^2 < 1\},$ $Y=(-1,1),$ and $f$ the
projection on the first factor.  Then there is no neighborhood of
$f^{-1}(0)$ on which $f$ becomes a trivial fibration, since the
inverse image of any $y\neq 0$ has finite genus and hence cannot
contain an embedded copy of $f^{-1}(0)$.  (To prove the latter fact,
notice that the intersection form on $H_1(S,\integers)$ has infinite
rank, and that this then has to be true for any manifold in which $S$
is embedded as an open subset.)  }
\end{ex}

Here is a more exotic example, which we will use later on as a
counterexample to an alternative version of our main theorem.

\begin{ex}[Exotic $\reals^4$'s]
\label{ex-exotic}
{\em Let $R$ be a manifold whose underlying topological space is
$\reals^4$, but which carries an exotic differentiable structure, and
let $h:R\to \reals$ be the length-squared function.  Although $h$
might not be differentiable, it is certainly continuous, so
$X=\{(x,r)|h(r)x^2 < 1\}$ is an open subset of $(-1,1)\times R$ and is
therefore a smooth manifold.  Let $f$ be the projection on the first
factor; it is locally trivial topologically but not in the smooth
sense.  In fact, if $R$ is chosen so that the fibres of $f$ form what
Gompf and Stipsicz \cite{go-st:kirby} call a radial family, their Theorem 9.4.10 
implies immediately that $f$ is not locally semitrivial, because $R$
cannot be embedded into any of the other fibres.
}
\end{ex}

To end this section, we show that Theorem \ref{thm-submersion} implies
our main theorem for the case of groupoids with trivial isotropy.

\begin{cor}
\label{cor-submersion}
Given a submersion $f:X\to Y$, let $G:X\arrows Y$ be the groupoid
which is the equivalence relation $X\times_Y X$.  For $y\in Y$, if
$\calo=f^{-1}(y)$ is of finite type, then there is a neighborhood
$\calu$ of $\calo$ in $X$  such that $G_{\calu}$ is isomorphic to
$\calo\times\calo\times f(\calu)\arrows \calo\times f(\calu)$, 
the product of the pair groupoid
$\calo\times\calo$ with the trivial groupoid $f(\calu)$.  This
product groupoid is in turn isomorphic to the action groupoid for the
action of $G_{\calo}$ on a neighborhood of the zero section in
the normal bundle $N\calo$.  
\end{cor}
\pf
The diffeomorphism $(\rho,f):\calu\to \calo\times f(\calu)$
of Theorem \ref{thm-submersion} gives an isomorphism from
$G_{\calu}$ to $(\calo\times f(\calu))\times_{f(\calu)}
(\calo\times f(\calu))$, which is isomorphic to the product groupoid
in the statement of the corollary. 

Furthermore, the derivative $Tf$ induces an isomorphism between $N\calo$ and
$\calo \times T_y Y$.  With respect to this isomorphism, the action of
$G\calo=\calo\times\calo$ on $N\calo$ is just given by its action on
$\calo$, with the trivial action on $T_y Y$, so that the action
groupoid is $\calo\times\calo\times T_y Y$. Finally, after shrinking
$\calu$ if necessary, we may identify $f(\calu)$ with a neighborhood
of zero in $T_y Y$, which gives the required isomorphism of
$G_{\calu}$ with the action groupoid.  
\nopagebreak
\qed

\section{{\'E}tale groupoids}
\label{sec-etale}
Recall that a Lie groupoid $G\arrows X$ 
is {\bf {\'e}tale} if its target (or source)
map is a local diffeomorphism.  The following theorem is the specialization
(slightly strengthened)
of Theorem \ref{thm-main} to the case where the groupoid is {\'e}tale and the
orbit consists of a single point.  (See also \cite{mo-mr:notes} and
\cite{mo-pr:orbifolds}.)

\begin{thm}
\label{thm-etale}
Let $G\arrows X$ be a proper {\'e}tale groupoid with fixed point $x\in X$.
Then there is a neighborhood $\calu$ of $x$ in $X$ such that the
restriction of $G$ to $\calu$ is isomorphic to the restriction of the
action groupoid $G_x \times T_x X \arrows T_x X$ to a neighborhood
of zero in $T_x X$.  If the groupoid is source-locally trivial, the
neighborhood $\calu$ can be taken to be $G$-invariant.
\end{thm}
\pf Since $G$ is proper and {\'e}tale, $G_x$ is a finite group.  Since $x$ is a
fixed point, $G_x$ is also the fibre $\alpha^{-1}(x).$ Hence, by (a
very simple case of) Theorem \ref{thm-submersion}, there is a
neighborhood $\calv$ of $G_x$ in $G$ such that the restriction of
$\alpha$ to $\calv$ is a trivial fibration.  $\calv$ is then a
disjoint union of finitely many open subsets of $G$ which are mapped
diffeomorphically by $\alpha$ to a neighborhood $\calu_1$ of $x$.  We
may assume that $\calu_1$ is connected.  

Let $r$ be the unique
continuous retraction from $\calv$ to $G_x$.  This map is a
``homomorphism'' (in quotes because $\calv$ is not necessarily a
subgroupoid of $G$) in the sense that $$r(g^{-1}h)=(r(g))^{-1}r(h))$$
for all $(g,h)$ in the fibre product $\calv \times_{\calu_1} \calv.$
To prove this, we note that each side of the displayed equation
is a continuous function of $(g,h)$ taking values in the
discrete space $G_x$, and that any pair $(g,h)$ in the fibre product
can be connected by a path to a pair in $G_x\times G_x$, where the
equation is obviously satisfied.

Via the action of $G$ on $X$ (see Example \ref{ex-action}),
the components of $\calv$ 
define an action of $G_x$ by local diffeomorphisms of $X$ fixing $x$.  
By the Bochner linearization theorem \cite{bo:compact} for actions of
compact groups, this action is equivalent in a neighborhood of $x$ to
the linearized action of $G_x$ on $T_xX$.  In particular, we can find
a disc $\calu$ about $x$ which is invariant under the action, and
hence is invariant for $G_{\calu_1}$ (though not necessarily for $G$
itself).  The restriction $G_{\calu}$ is then isomorphic to the action
groupoid for $G_x$ acting on a neighborhood of the origin in $T_x X$.

If $G$ is source-locally trivial, it follows from Theorem
\ref{thm-stability} that 
$\calu$ contains an invariant neighborhood of $x$;
the isomorphism with a transformation
groupoid still holds there.
\qed

\begin{rmk}
{\em We note for later use (in the proof of Theorem
\ref{thm-slice}) that the groupoid morphism
$r:G_{\calu}\to G_x$ associated to the isomorphism between $G_{\calu}$
and an action groupoid
is a covering morphism in the sense of 
\cite{br-da-ha:topological}; i.e. $(r,\beta)$ is a diffeomorphism
from $G_{\calu}$ to $G_x\times\calu.$ }
\end{rmk}

\begin{rmk}
{\em The groupoid $G_3$ in Example \ref{ex-twocharts} is proper and
{\'e}tale but is not
isomorphic to an action groupoid on
any invariant neighborhood of its fixed point.}
\end{rmk}

\begin{rmk}
{\em We conjecture that Theorem \ref{thm-etale} extends to the
non-{\'e}tale case.  See the discussion in
Section 4 of \cite{we:linearization}. }
\end{rmk}

\section{Deformation of proper groupoids}
\label{sec-deformation}
In this section, we prove that (target,source)-preserving deformations
of regular proper groupoids are trivial.  The result will be used to
extend the linearization theorem from {\'e}tale to effectively {\'e}tale
groupoids.  
As is usual in such deformation
problems, we use cohomology.

\begin{thm}
\label{thm-deformation}
Let $\{m_t\}$ be a smooth family of proper, regular, groupoid
structures on $G$ over $X$, defined for $t\in [0,1]$, having a fixed
map $(\alpha,\beta):G\to X\times X$ as anchor and a fixed
$\epsilon:X\to G$ as identity section.  Then there is a family
$\{A_t\}$ of diffeomorphisms of $G$ such that $A_0$ is the identity,
$A_t\circ\epsilon =\epsilon,$ 
and $A_t$ is a groupoid isomorphism from $(G,m_t)$ to
$(G,m_0)$.  
\end{thm}

\pf
Since the anchor is fixed, the submanifold $G^{(2)}\subseteq G\times
G$ of composable pairs is independent of
$t$.  For $(g,h)\in G^{(2)}$, we denote $m_t(g,h)$ by $g*_t h$.  
The derivative $Y_t$ of $m_t$ with respect to $t$ is a
mapping from $G^{(2)}$ to $TG$ which lifts $m_t$.  In fact, for
each composable pair $(g,h)$, $Y_t(g,h)$ lies in
the subspace of $T_{g*_t h}G$ which is annihilated by $(T\alpha,T\beta)$.
Differentiating the associativity law for the products $m_t$ with
respect to $t$, we obtain the identity
\begin{equation}
\label{eq-diffassociative}
Y_t(g,h)*_t k + Y_t(g*_t h,k) = g *_t Y_t(h,k) + Y_t(g,h*_t k),
\end{equation}
where the operation $*_t$ applied to 
a tangent vector and an element of
$G$ denotes the derivative of right [or left] translation
for the multiplication $m_t$; this derivative acts on
vectors tangent to the fibres of $\beta$ [or $\alpha$].

Equation (\ref{eq-diffassociative}) is actually a cocycle condition.
To see this, we right-translate the vectors $Y_t(g,h)$ back to the
identity section, i.e. we define $c_t(g,h)=Y_t(g,h)*_t(g*_t h)^{-1}.$
Since the values of $Y_t$ are tangent to the fibres of both $\alpha$
and $\beta$, the same is true of the values of $c_t$; i.e. $c_t(g,h)$
belongs to the fibre at $\alpha(g*_t h)=\alpha(g)$ of the isotropy
subalgebroid $\frakb$ of the Lie algebroid $\frakg$ of $G$.  As a
vector bundle, $\frakb$ is independent of $t$, and for each $t$ there
is an adjoint action of $G$ on $\frakb$ defined by  $g\bullet _t v = g *_t
v *_t g^{-1}$.  (See also 
Appendix B of \cite{ev-lu-we:transverse}).  

The infinitesimal associativity law (\ref{eq-diffassociative}) for $Y_t$
becomes the following identity for $c_t$:
$$c_t(g,h) + c_t(g*_t h, k) = g\bullet_t c_t(h,k) + c(g,h *_t k).$$
This identity says precisely that $c_t$ is a 2-cocycle on $G$
with values in the representation bundle $\frakb$.  

We now appeal to Proposition 1 of \cite{cr:differentiable}, which
establishes the triviality of the higher cohomology of any proper
groupoid with coefficients in any representation.  This provides us
with a family $\{b_t\}$ of 1-cochains whose coboundaries are $c_t$;
i.e. $b_t:G\to \frakb$ with $b_t(g)\in \frakb_{\alpha(g)}$ and 
\begin{equation}
\label{eq-coboundary}
b_t(g)+g\bullet_t b_t(h) - b_t(g*_t h) = c_t(g,h).
\end{equation}
Although it is not stated explicitly in \cite{cr:differentiable}, it
follows from the proof of Proposition 1 that the
``primitive'' $b_t$ can be chosen to depend smoothly on $t$.
(Alternatively, one can apply Proposition 1 to the single groupoid
$G\times [0,1]\arrows X\times [0,1]$
obtained by combining all the $(G,m_t)$.)  Furthermore, the
fact that all $m_t$ agree along the identity section implies that
$Y_t(g,g)=0$ whenever $g$ is an identity element, and then the
construction of $b_t$ shows that $b_t(g)=0$ as well.  

Reversing the translation procedure which led from $Y_t$ to $c_t$, we
construct from $b_t$ the vector field $X_t$ on $G$ defined by
$X_t(g)=b_t(g)*_t g$.  Since these vector fields vanish on the identities,
the family $\{X_t\}$ may be integrated at least locally to a smooth
family $\{A_t\}$ of diffeomorphisms of $G$ which fix the identities
and which commute with the anchor $(\alpha,\beta)$.  Since
the anchor is proper, the integration can be done globally.  Finally,
the coboundary relation (\ref{eq-coboundary}) is just the
differentiated version of the statement that $A_t$ is, for all $t$, a
groupoid homomorphism from $(G,m_t)$ to $(G,m_0)$.  
\qed
\begin{ex}
\label{ex-groups}
{\em
When $X$ is a point, Theorem
\ref{thm-deformation} implies the stability of group structures on
compact Lie groups and, with a ``multiparameter $t$,''
the local triviality of smooth
bundles of Lie groups.  In particular, the isotropy subgroupoid of a
regular proper groupoid is such a bundle, and hence all the isotropy groups
over a connected component of the base are isomorphic.  

Similarly,
if we have a locally trivial bundle of compact groups over a base manifold $Y$
acting freely on a locally trivial fibration $X\to Y$,
the action
groupoids form a locally trivial bundle of proper regular groupoids
over $Y$, and the orbit spaces
form a locally trivial bundle over $Y$.  
}
\end{ex}

\begin{rmk}
{\em 
It seems likely that Theorem \ref{thm-deformation} remains valid
even if the proper groupoid $G$ is not regular.  In this case, the
cohomology problem lives in a family of vector spaces which is not a
smooth bundle, but this should not be a serious difficulty.

Finally, we raise the question of deformability of proper groupoids,
regular or not, without the restriction that the anchor be fixed.  It
seems possible that a rigidity theorem like Theorem
\ref{thm-deformation} might still hold; it would be related to the
Reeb stability theorem.  
}
\end{rmk}

\section{Effectively {\'e}tale groupoids}
\label{sec-effectively}
A groupoid $G$ over $X$ is {\'e}tale when its maximal source-connected
subgroupoid is trivial.  We say that $G$ is {\bf effectively {\'e}tale}
when its maximal source-connected subgroupoid acts trivially on $X$;
equivalently, $G$ is effectively {\'e}tale when the anchor $\cala(G)\to
TX$ of its Lie algebroid is identically zero.  

If $G$ is effectively {\'e}tale,          then its maximal
source-connected subgroupoid  is a bundle of groups $B$ which
is a normal subgroupoid of $G$, and the quotient $G/B$ is {\'e}tale.  In
other words, an effectively {\'e}tale groupoid is an extension of an
{\'e}tale groupoid by a bundle of groups.  When the groupoid is proper,
the bundle of groups is locally trivial (see Example
\ref{ex-groups}).  
\begin{thm}
\label{thm-effectively}
Let $G\arrows X$ be a proper, effectively
{\'e}tale groupoid with fixed point $x\in X$.
Then there is a neighborhood $\calu$ of $x$ in $X$ such that the
restriction of $G$ to $\calu$ is isomorphic to the restriction of the
action groupoid $G_x \times T_x X \arrows T_x X$ to a neighborhood
of zero in $T_x X$.  If the groupoid is source-locally trivial, the
neighborhood $\calu$ can be taken to be $G$-invariant.
\end{thm}

\pf
Let $B$ the the maximal source-connected subgroupoid of $G$.  Since
$G$ is proper, so is $G/B$, and hence we may apply Theorem
\ref{thm-etale} to find a neighborhood $\calu$ of $X$ such that
$(G/B)_{\calu}$ can be identified with the restriction of the action
groupoid $(G/B)_x \times T_x X\arrows T_x$ to a neighborhood $\calv$
of $0$ in $T_x X$.  Note that $(G/B)_x\cong G_x/B_x$.  By choosing a
representative of each coset of $B_x$, we also obtain a
diffeomorphism (though generally not a group isomorphism) between $G_x$
and $B_x \times G_x/B_x$.  

By Theorem 5.1, $B$ is a locally trivial bundle of groups; its typical
fibre is $B_x$.  Hence, as a manifold, $G_\calu$ may be identified with 
$B_x \times G_x/ B_x\times \calv$, and the target and source maps
of its groupoid structure factor through those of the action groupoid
$G_x / B_x \times \calv \arrows \calv$.  Recalling that the groupoid
multiplication on the latter has the form $(g,hy)(h,y)=(gh,y)$, we
find that the multiplication on $B_x \times G_x/ B_x\times \calv$
coming from the groupoid structure on $G$ must have the form 
$m((a,g,hy),(b,h,y))=(M(a,b,g,h,y),gh,y)$, where $M:B_x \times B_x
\times G_x/B_x \times G_x/B_x\times \calu\to B_x$ is a smooth map.
When $y=0$ (corresponding to the fixed point $x$ in $\calu$), we have
$(M(a,b,g,h,0),gh)= (a,g)\cdot(b,h)$, the product in the isotropy group $G_x$.

We may now construct for $t\in [0,1]$
the smooth 1-parameter family of multiplications $m_t$
all having the anchor $(\alpha,\beta)(b,h,y)=(hy,y)$ by the
formula
$$m_t((a,g,hy),(b,h,y))=(M(a,b,g,h,ty),gh,y).$$  These multiplications
are all associative.  To see this, we note first that $t(hy)=h(ty)$ since the action of $G_x$
on $T_x X$ is linear, so 
$$m((a,g,thy),(b,h,ty))=(M(a,b,g,h,ty),gh,ty),$$
which shows that, for $t\neq 0,$ $m_t$ is the pullback of the
associative operation $m$ by the diffeomorphism
$(a,g,y)\mapsto (a,g,ty).$  On the other hand, for $t=0$, we have
$$m_0((a,g,hy),(b,h,y))=(M(a,b,g,h,0),gh,y)=((a,g)\cdot (b,h),y)$$
which is the product in the action groupoid $G_x\times\calv\arrows
\calv$.  We are thus in a position to apply Theorem \ref{thm-deformation}, 
which gives an isomorphism between this action
groupoid and the groupoid with multiplication $m_1$, which is just
$G_{\calu}$ itself.

If $G$ is assumed source-locally trivial, then, by Theorem
\ref{thm-stability}, $x$ is a stable fixed point from $G$, so the
neighborhood $\calu$ contains a $G$-invariant neighborhood on
which we still have an isomorphism with the linear approximation.
\qed

\section{The slice theorem}
\label{sec-slice}

A {\bf slice} through a point $x$ on an orbit $\calo$ of a groupoid
$G\arrows X$ will be defined simply as a  submanifold
$\Sigma$ of $X$ which meets $\calo$ only at $x$, with $T_x X = T_x \Sigma \oplus
T_x\calo.$   Only a small neighborhood of $x$ in a slice will be of
interest, and we can choose the neighborhood small enough so that it
is everywhere transverse to the orbits of $G$, so that the restriction
$G_{\Sigma}$ 
is again a Lie groupoid.  Since, by Proposition \ref{prop-closed}, 
$\calo$ is a closed submanifold, we can also suppose that
$\Sigma$ intersects $\calo$ only at $x$, so that $G_{\Sigma}$ has
$x$ as a fixed point.

In this section, we will show that a proper groupoid can be linearized
around the orbit $\calo$ if its restriction to a slice can be linearized around 
$x\in \calo$, and if $\calo$ is of finite type.  Combined with Theorem
\ref{thm-effectively}, this slice theorem will immediately imply our main
theorem.

\begin{thm}
\label{thm-slice}
Let $G\arrows X$ be a proper groupoid, and let $\calo$ be
an orbit of $G$ which is a manifold of finite type.  Suppose that the
restriction of $G$ to a slice through $x\in \calo$ is isomorphic
to the restriction of the action groupoid $G_x\times N_x\calo \arrows
N_x\calo$ to a neighborhood of zero.
Then there is a neighborhood $\calu$ of $\calo$ in
$X$ such that the restriction of $G$ to $\calu$ is isomorphic to the
restriction of the action groupoid $G_{\calo}\times_{\calo}
N\calo\arrows N\calo$ to a neighborhood of the zero section.
\end{thm}

\pf
The proof involves several steps, the basic idea being to apply
Theorem \ref{thm-submersion} to the restriction of $\alpha$ to
$G_{\Sigma X}$, where $\Sigma$ is a slice.

\noindent
{\bf Step 1.} 
Let $\Sigma$ be a slice as in the statement of the theorem, assumed small
enough so that it is everywhere transverse to the orbits of $G$
and intersects $\calo$ only at $x$, so that $x$ is a fixed point
of $G_{\Sigma}$.    $G_x$ is the isotropy group of $x$ in both
$G$ and $G_{\Sigma}$, and the natural identification of $N_x\calo$
with $T_x{\Sigma}$ is $G_x$-equivariant.

By assumption, the restriction $G_{\Sigma}$ is isomorphic, via a
retraction of groupoids $r:G_{\Sigma}\to G_x$ and an open embedding
$i:\Sigma\to T_x \Sigma$,
to the restriction of the action groupoid $G_x\times T_x \Sigma\arrows
T_x \Sigma$ to a neighborhood of zero.  We may assume that this neighborhood
is invariant for the action groupoid.  We then have an action of $G_x$
on $\Sigma$ as well as an embedding of $G_x$ into the group of bisections
of $G_{\Sigma}$.

\noindent {\bf Step 2.} Since the target map $\alpha$ is a submersion,
 $G_{\Sigma X}=\alpha^{-1}(\Sigma) $ is a closed submanifold of
$G_{\calu}$, where $\calu=G\Sigma$.  
The compact group $G_x$ acts on $\Sigma$ as mentioned in Step 1, and it also
acts freely on $G_{\Sigma X}$ by left translations via the embedding of $G_x$ into
the group of bisections of $G_{\Sigma}.$  
The restricted
submersion $\alpha :G_{\Sigma X}\to \Sigma$ is equivariant with respect to
these actions.  Applying Theorem \ref{thm-submersion}, we find a
$G_x$-equivariant local trivialization of $\alpha$ on a neighborhood
$\calv$ of $G_{x X}= \alpha^{-1}(x)$ in $G_{\Sigma X}$  That is,
there is a $G_x$-equivariant retraction $\rho:\calv\to G_{xX} $ such that
$(\rho,\alpha):\calv\to G_{xX} \times \Sigma$ is a diffeomorphism.  We will
denote the inverse of this diffeomorphism by $\Phi$.  

The $G_x$ orbits in $\calv$ are just the fibres of the source map
$\beta$, so if we let $\calu$ be $\beta (\calv)$, we obtain by
equivariance a retraction, also to be denoted by $\rho$, from $\calu$
to $\beta(G_{xX})=\calo$.  

Note that $\beta:\calv\to\calu$ and $\beta: G_{xX}\to \calo$ are
 principal bundles with structure group
$G_x$ (contrary to the usual conventions, the structure
group is acting on the left).  

\noindent
{\bf Step 3.}  In this step, we will construct a retraction of
 groupoids
 $R:G_{\calu} \to G_{\calo}$.

Let $p\in G_{\calu}$.  Since $\alpha(p)$ and $\beta(p)$ both belong to
$\calu=\beta (\calv)$, we can find $h$ and $k$ in $\calv$ such that
$\beta(h)=\beta(p)$ and $\beta(k)=\alpha(p)$.    The product
$kph^{-1}$ is then defined and, since $\alpha(kph^{-1})=\alpha (k)$ and
$\beta(kph^{-1})=\beta(h^{-1}) =\alpha(h)$, $kph^{-1}$ lies in
$G_{\Sigma}$; hence, we have the element $r(kph^{-1})$ of 
$G_x$.
Now we define $R(p)$ to be $\rho(k)^{-1} r(kph^{-1}) \rho(h)$.  The target and
source of $R(p)$ are then
$\alpha(R(p))=\beta(\rho(k))=\rho(\beta(k))=\rho(\alpha(p))$
 and
$\beta(R(p))=\beta(\rho(h))=\rho(\beta(h))=\rho(\beta(p))$
respectively; in particular, $R(p)$ belongs to $G_{\calo}$. 

We must show that $R(p)$ is independent of the choices we made.  If
$h$ and $k$ are replaced by $h'$ and $k'$, we have $h'=bh$ and $k'=ck$
for $b$ and $c$ in $G_{\Sigma}$.  
Now we have
$k'ph'^{-1}=c(kph^{-1})b^{-1}$, 
so $r(k'ph'^{-1})=r(c)r(kph^{-1})r(b)^{-1}$, so 
$\rho(k')^{-1} r(k'ph'^{-1}) \rho(h') = \rho(ck)^{-1}
r(c)r(kph^{-1})r(b)^{-1}
\rho(bh).$  By
the equivariance of $\rho$, this becomes 
$$\rho(k)^{-1}r(c)^{-1}
r(c)r(kph^{-1})r(b)^{-1}
r(b)\rho(h)=\rho(k)^{-1}r(kph^{-1})\rho(h)$$
as before, so $R(p)$ is well
defined.  

To see that $R$ is a retraction, if we assume that $p$ belongs to
$G_{\calo}$, we may choose $h$ and $k$ in $G_{xX}$, the image of the
retraction 
$\rho$, so that $kph^{-1}\in G_X,$ and hence
$R(p)=\rho(k)^{-1} r(kph^{-1}) \rho(h)
=k^{-1}kph^{-1}h=p.$

Finally, we will show that $R$ is a groupoid homomorphism, with the
map on objects being $\rho:\calu\to\calo.$  We have already seen that
$R$ and $\rho$ are compatible with the source and target maps.  If
$p$ and $q$ are in $\calu$ with
$\beta(p)=\alpha(q),$ and we choose $h$ and $k$ as above to compute
$R(p)$, then we may use some  $g$ and the same $h$ to compute $R(q)$, while
$g$ and $k$ may be used to compute $R(pq)$.  With
$a=kph^{-1} $ and $b=hqg^{-1}$, we find $ab=kpqg^{-1}$, and so
$R(p)R(q)=\rho(k)^{-1} r(a)\rho(h) \rho(h)^{-1}r(b) \rho(g)=\rho(k)^{-1}
r(ab) \rho(g) =R(pq)$.

\noindent
{\bf Step 4.}  To identify $G_{\calu}$ with an action groupoid, we
will show that the homomorphism $R$ is a covering morphism in the
sense of \cite{br-da-ha:topological}.  This means that we must show
that the map $(R,\beta)$ is a diffeomorphism from $G_{\calu}$ to 
$G_{\calo}\times_{\calo}\calu,$ the fibre product with respect to the
pair $(\beta,\rho)$.  To do so, we will construct an inverse map $\Psi$.

Let $(m,z)\in G_{\calo}\times_{\calo}\calu.$  Since
$\calu=\beta(\calv)$, we may write $z=\beta(h)$ for some $h\in \calv$;
in particular, $\alpha(h)\in \Sigma.$   Then
$\beta(\rho(h))=R(\beta(h))=R(z)=\beta(m),$ since $(m,z)$ is in the
fibre product, so the product $u=\rho(h)m^{-1}$ is defined, with
$\alpha(u)=x$ and
$\beta(u)=\alpha(m)\in\calo.$  Recalling that $(\rho,\alpha):\calv\to
G_{xX} \times \Sigma$ is a diffeomorphism with (equivariant)
inverse $\Phi,$ we set
$\Psi(m,z)=\Phi(u,\alpha(h))^{-1}h,$ i.e.
$ \Psi(m,z)=\Phi(\rho(h)m^{-1},\alpha(h))^{-1}h.$

To see that $\Psi$ is well defined, we replace $h$ by $h'=ah$, where
$a\in G_{\Sigma}$.  Carrying out the construction of the previous paragraph
with this new choice, we have 
$u'=\rho(ah)m^{-1}=r(a)\rho(h)m^{-1}$ and 
$$\Phi(u',\alpha(h'))=\Phi(r(a)\rho(h)m^{-1},\alpha(ah))=
\Phi(r(a)\rho(h)m^{-1},a\alpha(h))$$
$$=a\Phi(\rho(h)m^{-1},\alpha(h))=a\Phi(u,\alpha(h)).$$
Thus $\Phi(u',\alpha(h'))^{-1}h'=\Phi(u,\alpha(h))^{-1}a^{-1}ah
=\Phi(u,\alpha(h))^{-1}h,$ so $\Psi$ is well-defined.  

Now we show that $\Psi$ is indeed an inverse to $(R,\beta)$.  For $p\in
G_{\calu},$ we choose $h$ and $k$ as in Step 3 and set $a=kph^{-1}$; then 
$R(p)=\rho(k)^{-1} r(a) \rho(h).$  Now, since $\beta(h)=\beta(p)$,
$$\Psi(R(p),\beta(p))=\Phi(\rho(h)R(p)^{-1},\alpha(h))^{-1}h$$
$$=(\Phi(\rho(h)\rho(h)^{-1}r(a)^{-1}\rho(k),\alpha(h))^{-1}h$$
$$=(\Phi(r(a)^{-1}\rho(k),\alpha(h))^{-1}h
(\Phi(r(a)^{-1}\rho(k),a^{-1}a\alpha(h))^{-1}h$$
$$=(\Phi(\rho(k),\alpha(ah))^{-1}ah=
(\Phi(\rho(k),\alpha(kp))^{-1}kp=k^{-1}kp=p.$$

In the other direction, for $(m,z)\in G_{\calo}\times_{\calo}\calu,$
with  $z=\beta(h),$ we have 
$$\beta(R(\Psi(m,z))=\beta(\Phi(\rho(h)m^{-1},\alpha(h))^{-1}h)=\beta(h)=z,$$ 
while
$$R(\Psi(m,z))=R(\Phi(\rho(h)m^{-1},\alpha(h))^{-1}h).
$$

To show that the last expression is equal to $m$, thus completing 
the proof, we may choose for the $h$ in the
definition of $R$ the $h$ which we used to define $\Psi(m,z).$
 For the $k$ in
the definition we take
$\Phi(u,\alpha(h))=\Phi(\rho(h)m^{-1},\alpha(h)).$  Then 
$$R(\Psi(m,z))
=\rho(k)^{-1}r(k\Psi(m,z)h^{-1})\rho(h)$$
$$=\rho\Phi(u,\alpha(h))^{-1} 
r(\Phi(u,\alpha(h))\Psi(m,z)h^{-1})\rho(h)$$
$$=u^{-1}r(hh^{-1})\rho(h)=u^{-1}\rho(h)=m.$$
\qed

\section{Proof of the main theorem}
\label{sec-mainproof}

For convenience, we restate the main theorem.

\bigskip
\noindent
{\bf Theorem 4.1}
{\em
Let $G\arrows X$ be a regular, proper Lie groupoid, and let $\calo$ be
an orbit of $G$ which is a manifold of finite type.  Then 
there is a neighborhood $\calu$ of $\calo$ in
$X$ such that the restriction of $G$ to $\calu$ is isomorphic to the
restriction of the action groupoid $G_{\calo}\times_{\calo}
N\calo\arrows N\calo$ to a neighborhood of the zero section in $N\calo$. 
}

\bigskip
\pf

Let $\Sigma$ be a slice through $x\in \calo.$   If $\Sigma$ is
chosen small enough, then $x$ is a fixed point of $G_{\Sigma}$, so
the anchor of the Lie algebroid $\cala(G_{\Sigma})$ is zero at
$x$.  Since $G$ is regular, so is $G_{\Sigma}$, and hence the anchor of 
$\cala(G_{\Sigma})$ is identically zero; i.e. $G_{\Sigma}$ is
effectively {\'e}tale.  By Theorem \ref{thm-effectively}, we can choose
$\Sigma$ small enough so that $G_{\Sigma}$ is locally isomorphic
to its linearization at $x$.  By Theorem \ref{thm-slice}, $G$
is isomorphic to its linearization along $\calo$.
\qed
We close with some remarks.

\begin{ex}
\label{ex-exoticcross}
{\em 

It was tempting to substitute an assumption of
source-local triviality of $G\arrows
X$ for the hypothesis in the main theorem that $\calo$ be of
finite type.  The following example shows that this is not
possible.

We begin with the submersion $f:X\to (-1,1)$ of Example 
\ref{ex-exotic}.  Let $G\arrows X$ be the groupoid which is the
product of the equivalence relation $X\times_{(-1,1)} X \arrows X$
and the group $SU(2)$, where the latter is considered as a
groupoid over a one-point base.  The fibre over $(x,r)$
of the source map of $G$ is $f^{-1}(r)\times SU(2)$, which is
topologically the product $\reals^4 \times SU(2)$, and the source map is 
locally trivial as a topological fibration.   By Theorem 2 of
\cite{me:submersions1}, the source map is locally trivial in the
differentiable sense as well; i.e. $G\arrows X$ is source-locally
trivial.  

On the other hand, $G\arrows X$ is not isomorphic to its
linearization around $\calo=f^{-1}(0)$ in any neighborhood of
$\calo$, since if it were, the restriction of $f$ to such a
neighborhood would be locally trivial, and we saw in Example
\ref{ex-exotic} that this is not the case.
}
\end{ex}

Finally, we would also like to mention an alternative approach to
understanding proper regular groupoids, due in part to I. Moerdijk (private
communication).  Instead of first
restricting $G\arrows X$ to a slice to get an effectively \'etale groupoid and
then dividing by the identity component of the isotropy, we may
first divide $G$ itself by the identity component $C$ of its
entire isotropy, which turns out to be a smooth bundle of compact
groups.  The quotient $G/C$ is then a {\bf foliation groupoid}
over $X$, i.e. it is a groupoid for which the Lie algebroid anchor
is injective.  (Note that the orbits, or ``leaves,'' of $G/C\arrows X$ are not 
necessarily connected.)
The foliation groupoid can be analyzed via the
slice theorem in terms of its restriction to a slice, which is the
same \'etale groupoid as was obtained in the first approach.
Finally, the original groupoid $G\arrows X$ may be seen as an
extension of a foliation groupoid by a bundle of compact groups.
The bundle is locally trivial (Example \ref{ex-groups}), and the
extension is then classified by a degree 2 cohomology class of the 
foliation groupoid with values in the bundle of groups.   If we
restrict to a neighborhood of an orbit, we can again use a
deformation argument, somewhat more complicated than that in
Theorem \ref{thm-effectively}, to recover the main linearization
theorem.

\appendix

\section{Appendix: Proper mappings}
\label{sec-propermapping}
Perhaps the most common definition of properness for mappings between
topological spaces is:
\begin{dfn}
\label{dfn-propercompact}
A mapping $f:X\to Y$ between topological spaces is proper if it is
continuous and if $f^{-1}(A)$ is compact in $X$ for every compact
subset $A$ of $Y$.
\end{dfn}
There are two other definitions of properness which are equivalent
to this one when $X$ and $Y$ are Hausdorff spaces, but which differ in
general.  Since many interesting groupoids, such as holonomy groupoids
of foliations, may not be Hausdorff, we mention these other definitions. 

James \cite{ja:general} defines properness of a map in the following
way.

\begin{dfn}
\label{def-properproduct}
A mapping $f:X\to Y$ between topological spaces is proper if it is
continuous and if, for every topological space $Z$, the product
mapping $f\times 1_Z:X\times Z \to Y \times Z$ is closed, in the sense
that it maps closed sets to closed sets.
\end{dfn}

Actually, James calls such mappings ``compact,'' but we will not use
this term.  James actually goes on to define a
topological space $X$ to be compact if the map from $X$ to a point is
proper, and then he proves that this definition is equivalent to the
usual one in terms of open coverings.  He also proves:

\begin{prop}
\label{prop-compactfibres}
A continuous mapping is proper if and only if it is closed and the
inverse image of each point is compact.
\end{prop}

Another definition of properness is given by Crainic and Moerdijk in
\cite{cr-mo:homology}. 

\begin{dfn}
\label{def-properfibre}
A continuous mapping $f:X\to Y$ is proper if: (i) the image of the
diagonal $X\to X\times_Y X$ is closed; and (ii) $f^{-1}(A)$ is compact
whenever $A$ is a compact subset of a Hausdorff open subset of $Y$.
\end{dfn}

This definition is adapted to the study of non-Hausdorff groupoids and
their associated operator algebras.  It leads easily to the conclusion that, if 
$G\arrows X$ is a proper groupoid, and $X$ is Hausdorff, then $G$ is Hausdorff
as well.  It is not clear how this definition relates to
the one which is expressed in terms of products.

For the reader's convenience, we also include here a proof of a standard
fact
used in Section \ref{sec-submersions}.

\begin{lemma}
\label{lemma-complete}
If $K$ is a compact group, any $K$-manifold $M$ admits a complete,
invariant riemannian metric.
\end{lemma}
\pf
It is clearly sufficient to prove the lemma under the assumption that
$M$ is connected.  This implies that $M$ admits a partition of unity
by compactly supported functions $M\to [0,1]$ which can be enumerated
$\pi_1, \pi_2 \ldots .$    The sum $\sum_{n=1}^{\infty} n\pi_n$ (or
the corresponding finite sum if $M$ is compact) is then a proper,
nonnegative function on $M$.  By averaging with respect to $K$ we
obtain an invariant function $\lambda$.  Using the compactness of $K$,
one shows easily that the averaged function is again proper. 

Now let $\langle~,~\rangle_0$ be any invariant riemannian metric on
$M$, and define the new metric $\langle~,~\rangle$ to be $(1+\langle \nabla_0
\lambda, \nabla_0 \lambda\rangle)\langle~,~\rangle_0$.  Then the gradient
$\nabla \lambda$ of $\lambda$ with respect to $\langle~,~\rangle$ has
length everywhere less than $1$.  It follows that $\langle~,~\rangle$
is complete, since the proper function $\lambda$ is bounded on any
curve of finite length, which implies that closed bounded subsets of
$M$ are compact.
\qed

\section{Appendix: Manifolds of finite type}
\label{sec-finite}

There does not seem to
be a standard name for the following concept.

\begin{dfn}
\label{dfn-finite}
Let $K$ be a compact group.  A $K$-{\bf manifold of finite type} is a
$K$-manifold $M$ which admits a proper $K$-invariant function whose
critical points form a compact set.  When $K$ is the trivial group, we
simply say that $M$ is a manifold of finite type.
\end{dfn}

Here are some elementary observations about manifolds of finite type.

By squaring a given proper function, we can arrange that the proper
function in the definition above take values in $[0,\infty)$.  If $K$
is finite, we can arrange that the function have finitely many
critical points, and even that they be nondegenerate if $K$ is
trivial.  We also note that $M$ is a $K$-manifold of finite type if
and only if it is equivariantly diffeomorphic to the interior of a
compact $K$-manifold with boundary.

If $M$ is of finite type, then so is any $K$-equivariant bundle over
$M$ with compact fibres, e.g. a finite covering.  If $K$ is a compact
group acting freely on $M$, then $M$ is a $K$-manifold of finite type
if and only if $M/K$ is a manifold of finite type.  But note the
following example.

\begin{ex}
\label{ex-exotic1}
{\em There exist many manifolds which are homeomorphic
to $\reals^4$ and which are not of finite type.  (See
\cite{go-st:kirby}.  In fact, Gompf and Stipsicz remark on page 366
that the existence of an exotic $\reals^4$ of finite type would lead
to a counterexample to the differentiable Poincar\'e conjecture in
either dimension 3 or dimension 4.)
Let $Q$ be such a manifold of infinite type, and let 
$M=Q\times SU(2)$.  Then $M$ is simply connected
  and is simply connected at infinity, so it is of finite type, according
to Siebenmann \cite{si:obstruction}.  Let $SU(2)$ act on $M$
  by left translation on the second factor.  Then $M/SU(2)$ is $Q$,
  which is not of finite type. 
}
\end{ex}

\begin{rmk}
\label{rmk-coverings}
{\em  It would be interesting to know whether a finite group
can act freely on a manifold of
   finite type such that the quotient manifold is of infinite type.
An example might be used to extend Example \ref{ex-exoticcross} to
cover the \'etale case.
}
\end{rmk}

\end{document}